\documentclass[aip, amsmath,amssymb,reprint]{revtex4-2}
\usepackage{graphicx}
\usepackage{dcolumn}
\usepackage{bm}
\usepackage[utf8]{inputenc}
\usepackage[T1]{fontenc}
\usepackage{mathptmx}
\usepackage{verbatim}
\usepackage{amsfonts,amsthm,amsmath,amssymb,amscd,graphicx,color}

\usepackage{mathtools}
\usepackage{hyperref}
\usepackage{MnSymbol,wasysym}
\usepackage{mathtools}
\usepackage[toc,page]{appendix}
\usepackage{soul}

\makeatletter
\DeclareTextCommand{\textprime}{\encodingdefault}{%
  \mbox{$\m@th'\kern-\scriptspace$}%
}
\makeatother

\begin{document}

\preprint{AIP/123-QED}

\title[Measuring chaos in the Lorenz and R\"ossler models: Fidelity tests for reservoir computing.]
      {Measuring chaos in the Lorenz and R\"ossler models:\\ Fidelity tests for reservoir computing.}

\author{James Scully}
\affiliation{Neuroscience Institute, Georgia State University,
100 Piedmont Ave., Atlanta, GA 30303, USA.}
\author{Alexander Neiman}
\affiliation{Department of Physics and Astronomy, Ohio University, Athens, OH 45701, USA.}
\author{Andrey Shilnikov}
\affiliation{Neuroscience Institute and Department of Mathematics \& Statistics, Georgia State University,
100 Piedmont Ave., Atlanta, GA 30303, USA.}

\date{\today}

\begin{abstract}
This study is focused on the qualitative and quantitative characterization of chaotic systems  with the use of symbolic description. We consider two famous systems: Lorenz and R\"ossler models with their iconic attractors, and demonstrate that with adequately chosen symbolic partition three measures of complexity, such as the Shannon source entropy, the Lempel-Ziv complexity and the Markov transition matrix, work remarkably well for characterizing the degree of {\em chaoticity}, and precise detecting stability windows in the parameter space. \\
The second message of this study is to showcase the utility of symbolic dynamics with the introduction of a fidelity test for reservoir computing for  simulating the properties of the chaos in both models' replicas. 
The results of these measures are validated by the comparison approach based on one-dimensional return maps and the complexity measures. \end{abstract}
\maketitle

\textbf{ We employ the methods of qualitative theory and symbolic dynamics to measure  chaos and detect stability islands in one-parametric sweeps in the Lorenz and R\"ossler models, as well as  to compare chaotic properties in their reservoir computed surrogates. We seek to test that reservoir computing algorithms are able to pass all tests prepared for them: quantitive ones relaying on all three measures extracted from binary sequences such as block entropy, Lempel-Ziv complexity and Markov matrix structures, as well as qualitative ones based on return maps.  We hope that our algorithms and findings will be helpful for a broad interdisciplinary audience including specialists and beginners in dynamical systems and machine learning.}

\textbf{We dedicate this article to the memory of our colleague, teacher, and dear friend -- Professor Vadim Anischenko who made a fundamental contribution to the field of nonlinear dynamics. Among so many aspects of nonlinear systems, statistical properties of complex dynamics and methods of their characterization were of his continuous interest.}


\section{Introduction}

Finding effective characterization of complex time series is a pivotal task for the understanding of their underlying dynamics \cite{abarbanel1993analysis,kantz2004nonlinear,bradley2015nonlinear}.
The Lorenz and R\"ossler models are the classic examples of two types of deterministic chaos observable in various low-dimensional dynamical systems \cite{anishchenko1990complex}, respectively, with Lorenz-like attractors and spiral ones due to the Shilnikov saddle-focus. 
A qualitatively different mechanism of formation and structure of chaos are reflected in  distinct statistical properties of these systems \cite{anishchenko2003correlation,anishchenko2005statistical,anishchenko2001effect,anishchenko2004autocorrelation}. As such, the Lorenz and R\"ossler models serve as test-benches for testing and development of new tools in the field of nonlinear science.

The first goal of this paper is to showcase how one may measure the degrees of chaotic, homoclinic dynamics in such systems based upon the symbolic description and how well the proposed approach  complimentary agrees with the conventional one employing the Lyapunov exponents. We also argue that the symbolic approaches work exceptionally well to detect the stability windows in parametric sweeps of such systems. The second goal of our paper is to experiment to what degree machine learning tools such as reservoir computing may well learn to reflect qualitatively and quantitatively on the dynamical and probabilistic properties of original systems and surrogated ones.          

While neither model needs to be introduced to the nonlinear community, nevertheless let us first of all 
describe some of the key dynamical feathers of both systems, as well as how their chaotic dynamics can be translated into the symbolic description to generate long binary sequences to be further quantified and analyzed using a simple technique of partitioning the phase space or phase variables of these classic models.  
 
 The paper is organized as follows: first we discuss the Lorenz model, followed by the R\"ossler model, and introduce the suitable symbolic description for both. Next, we introduce the complexity measures using on the binary framework; they include block-entropies and source entropy of symbolic sequences, Lempel-Ziv (LZ) complexity. We argue that the entries of the Markov transition matrix can effectively indicate structurally stable dynamics in such pseudo-hyperbolic system, and therefore to detect the stability windows in the parameter space.  Next we compare these complexity measures with the largest Lyapunov exponent as effective computational indicators of chaos and periodic dynamics. Finally, we analyze and quantify the closeness of the original chaotic dynamics occurring in the Lorenz and R\"ossler model and their clones generated by various recurrent neural networks based on reservoir computing principles.      

\section{Models, symbolic partitions and binary framework}
 
\subsection{Lorenz model}

\begin{figure}[t!]
	\includegraphics[width=0.85\linewidth]{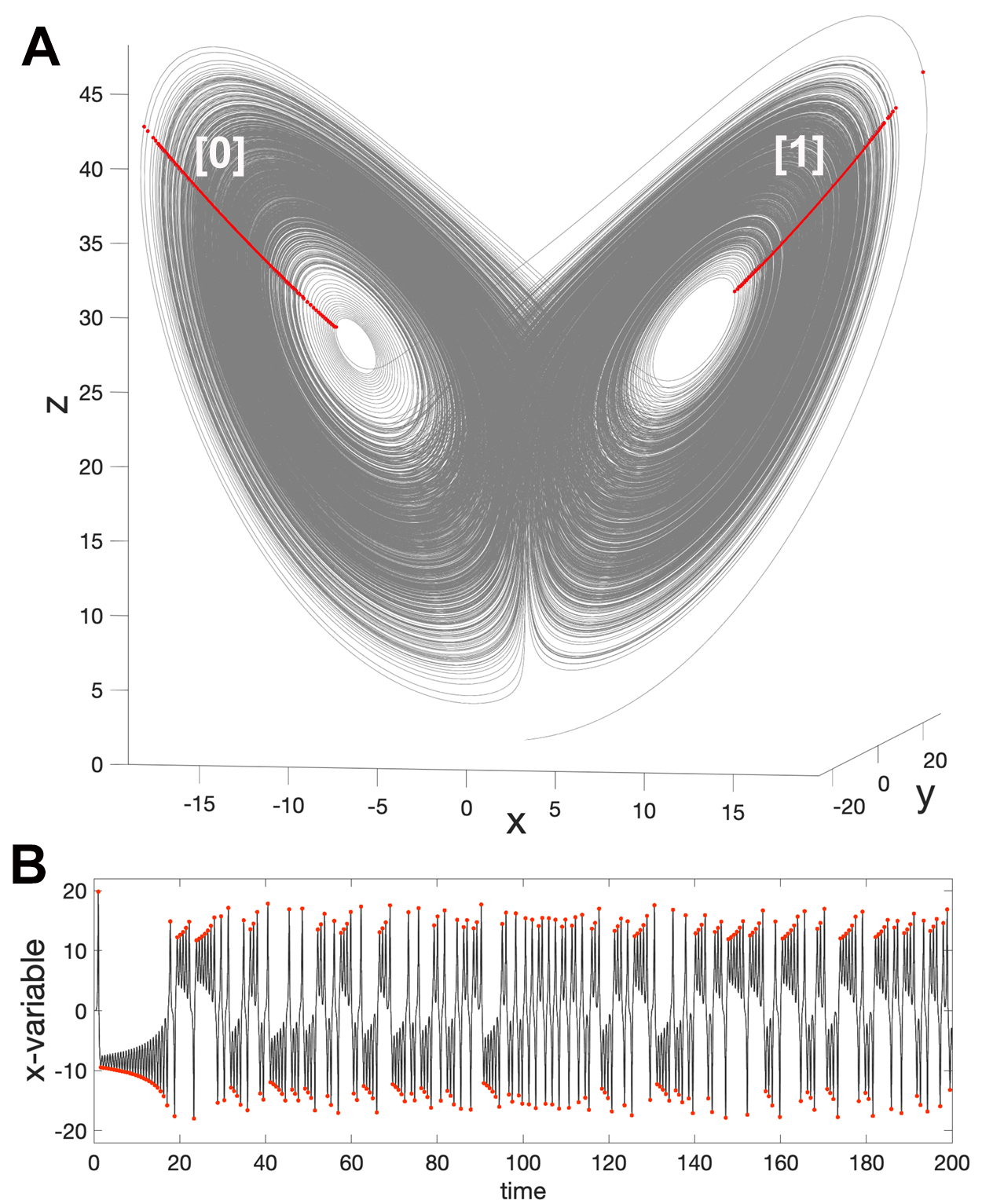}
	\caption{(A) The Lorenz attractor at the classic value $r=29$. The superimposed red dots defined by the $x$-variable critical events are well-aligned on some straight-line intervals transverse to the wings of the Lorenz butterfly in the phase space. (B) The $x$-variable plotted against time.  Local maxima and minima marked with red dots are detected to convert the $x$-dynamics into binary sequences using the simple rule: $\{\dot x =0\ |\, x>0\} \to ``0"$  and $\{\dot x=0\ |\, x<0\} \to ``1"$ in this case.}
	\label{fig1}
\end{figure}

\begin{figure}[t!]
	\includegraphics[width=0.99\linewidth]{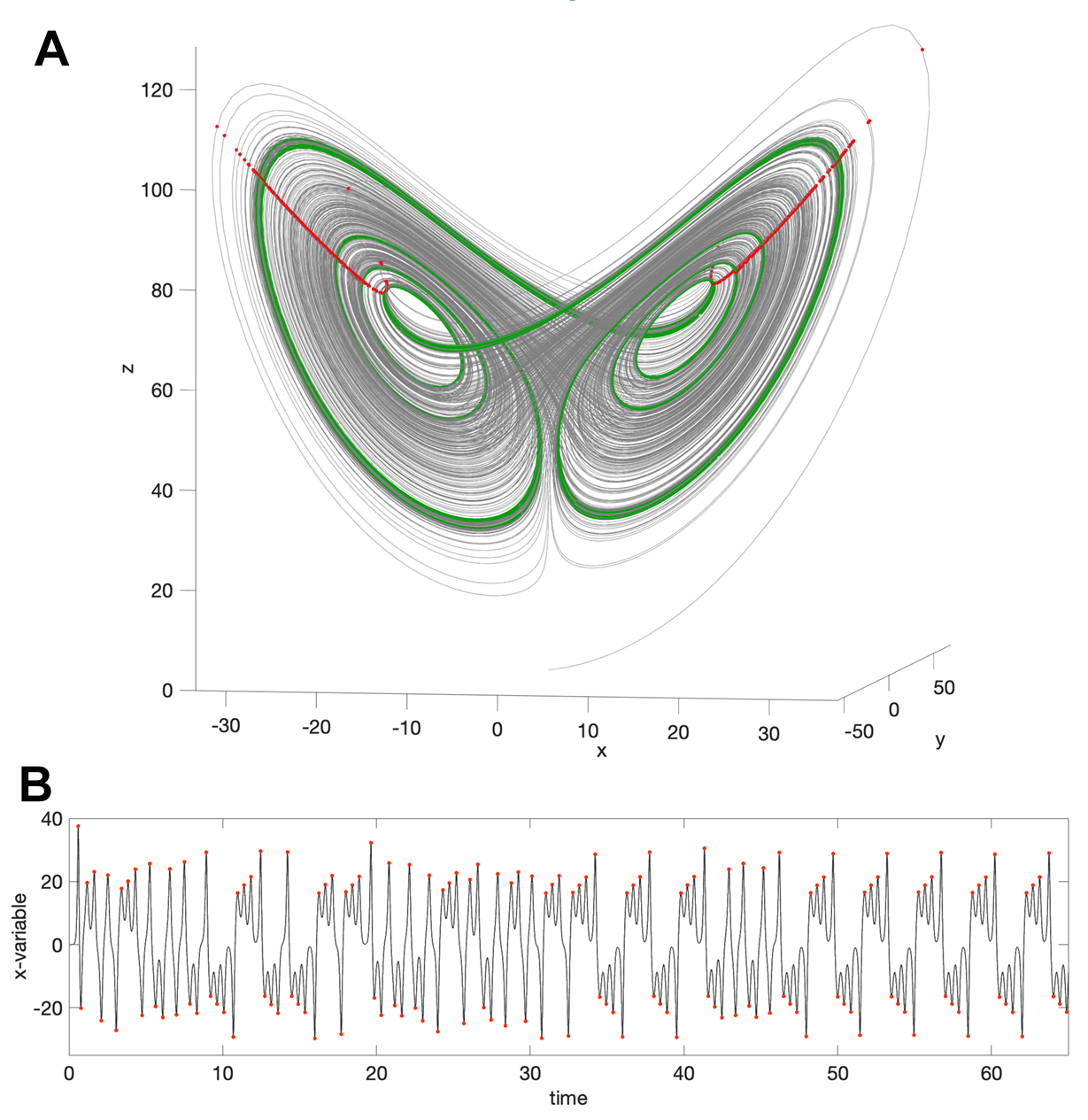}
	\caption{(A) Chaotic transient converging to a stable attractor (green), encoded as [$\overline{00001111}$], in the 3D phase space of the Lorenz model near a stability window at  $r\simeq$69.67. (B) The $x$-variable time traces passes through a periodic pattern with the Markov transition probabilities $p_{00}=3/4$ and $p_{01}=1/4$. }
	\label{fig2}
\end{figure}

The Lorenz equation or model is given by:
\begin{equation}
\dot x = -\sigma(x- y), \quad \dot y = r\,x -y +x\,z, \quad  \dot z = bx +x\,y,
\label{eq:Lorenz}
\end{equation}
 with $x,y$, $z$ being the phase variables, and $\sigma,\ r,\ b>0$ being bifurcation parameters; we will keep $\sigma=10$ and $b=8/3$ fixed through this study. For $r>28$ the model starts exhibiting chaotic behavior associated with an iconic butterfly-shaped strange attractor depicted in Figs.~\ref{fig1}-\ref{fig3} below. As this model is $\mathbb{Z}_2$-symmetric, i.e., it supports the group symmetry $(x,y,z) \leftrightarrow (-x,-y,z)$. This is manifested in the shape the Lorenz attractor shown  in the projection in Fig.~\ref{fig1}A, which is filled in with flip-flopping patterns of a single solution of  Eqs.~(\ref{eq:Lorenz}). This is well seen from Fig.~\ref{fig1}B representing a typical evolution of the $x$-variable in time. One can see from this  and similar traces shown in Figs.~\ref{fig2}B and \ref{fig3}B that switching $x$-patterns changes with variations of the bifurcation $r$-parameter. Specifically, one can see that the $x$-patterns in the last two figures become periodic after some chaotic transients. These correspond to the stable periodic orbits, shown in green, embedded in the chaotic attractors in the phase space as shown in Figs.~\ref{fig2}B and \ref{fig3}B.  On the other hand, for other $r$-parameters values, the Lorenz attractor remains chaotic as shown in Figs.~\ref{fig:L28}A-\ref{fig:L75}A and can be seen from time-progressions of the $z$-variable in Figs.~\ref{fig:L28}B-\ref{fig:L75}B. This property makes the Lorenz attractor {\em pseudo-hyperbolic} or a {\em quasi-attractor} \cite{ABS:1977,AfrShil1983,ABS:1983,Shilnikov1986,BykovShilnikov92,GTS93}. Without going into details, the first means that it constantly changes due to homoclinic bifurcations of its key contributor -- the saddle at the origin with two outgoing separatrices that fill in the butterfly of the Lorenz attractor in the given phase space projections. In contrast, a pseudo-hyperbolic attractor becomes a quasi-attractor with homoclinic tangencies causing the emergence of stable periodic orbits within it that remains such on $r$-parameter intervals, known as the stability windows or islands.       
 
  In this study, we will demonstrate how time-progressions of both variables can be used to compare and contrast the chaotic and periodic dynamics generated by the Lorenz model and to detect stability islands in parameter space.     

Let us first discuss the way symbolic, binary representations can be introduced to describe the flip-flop dynamics of the Lorenz model. Figure~\ref{fig1} illustrates the concept: every turn of a phase trajectory around the right "eye" in the wing of the butterfly, which is given by the condition $x'=0|~~ x>0$,  generates symbol "1" in the binary sequence. Otherwise, a turn around the left eye, when $x'=0|~~x<0$, adds "0" symbol to the  sequence. 
For example, the sequence $\{10000..10111...\}$ for the $x$-variable progression shown in Fig.\ref{fig1}B that corresponds to the right separatrix of the saddle at the origin in the 3D phase space. 
Such a sequence can be aperiodic/chaotic for the canonical parameter value $r=28$, used in all textbook on nonlinear dynamics. For other values, the sequence can become periodic with a repetitive block, e.g. [$\overline{00001111}$] of period 8, after some short or long transient, as in the case depicted in Fig.~\ref{fig2}. The periodic orbit in Fig.~\ref{fig3} has a shorter periodic block [$\overline{000111}$] of period 6, and so forth. 

This outlines the method and ultimate goal of the symbolic approach: first, one picks some typical, long steady-state trajectory of the model, 
and second, examines the binary code/sequence extracted from the $x$-variable progression. The question is how to determine efficiently whether the sequence is periodic or aperiodic, i.e., chaotic and what is the degree of chaoticity/complexity? 

 \begin{figure}[t!]
	\includegraphics[width=0.999\linewidth]{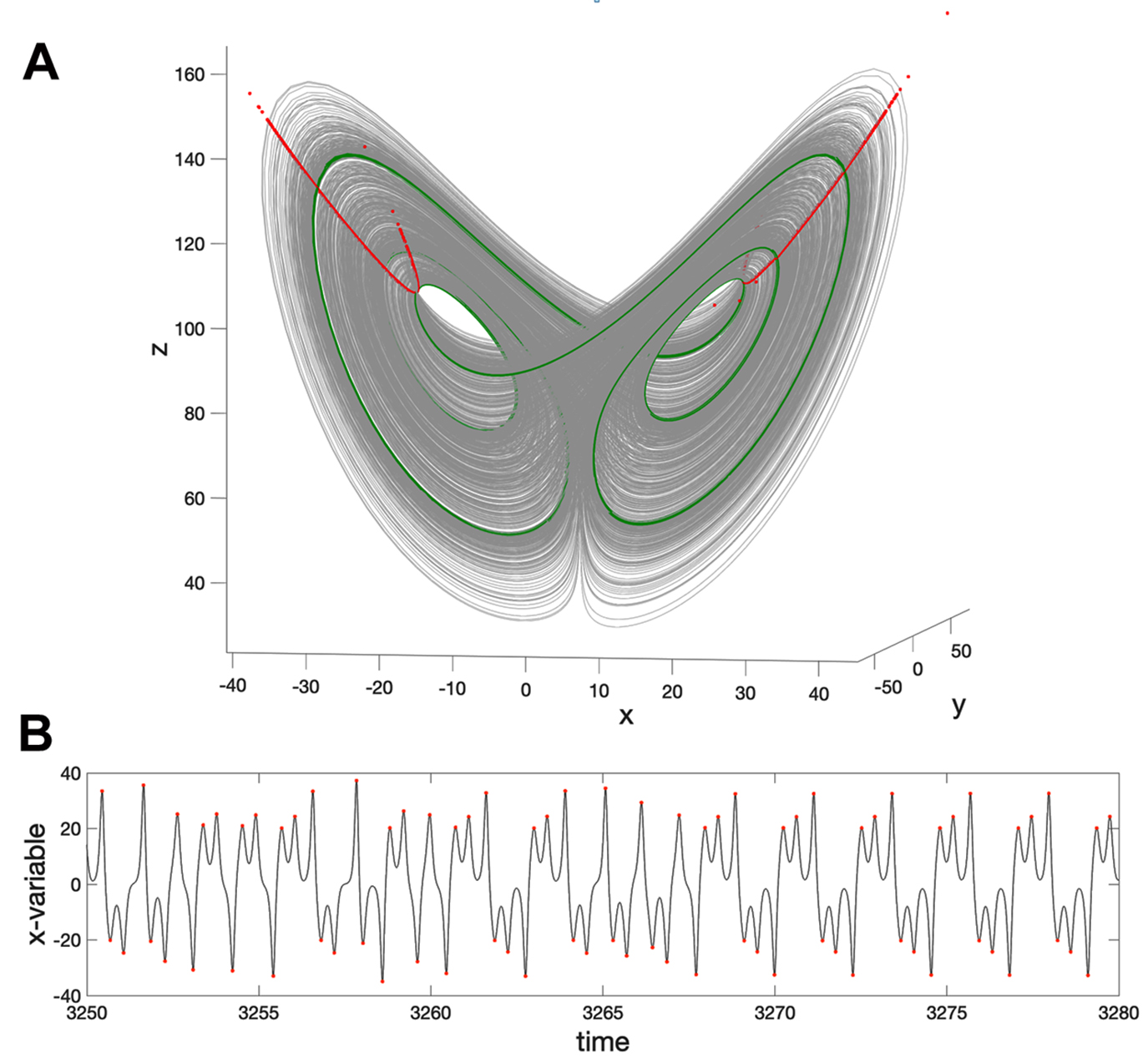}
	\caption{(A) Convergence to the stable periodic orbit (green), encoded as [$\overline{000111}$], after a long chaotic transient in the 3D phase space of the Lorenz model at  $r=92.5$. Superimposed red dots defined as critical events $\{x'=0|\, x>0$ and $x<0\}$ fill out two hooks on the bending wings of the butterfly in the phase space. (B) The $x$-variable plotted against time reveals the attracting periodic pattern with the Markov transition probabilities $p_{11}=2/3$ and $p_{10}=1/3$. }
	\label{fig3}
\end{figure}

\subsection{R\"ossler model}

 \begin{figure}[t!]
	\includegraphics[width=0.999\linewidth]{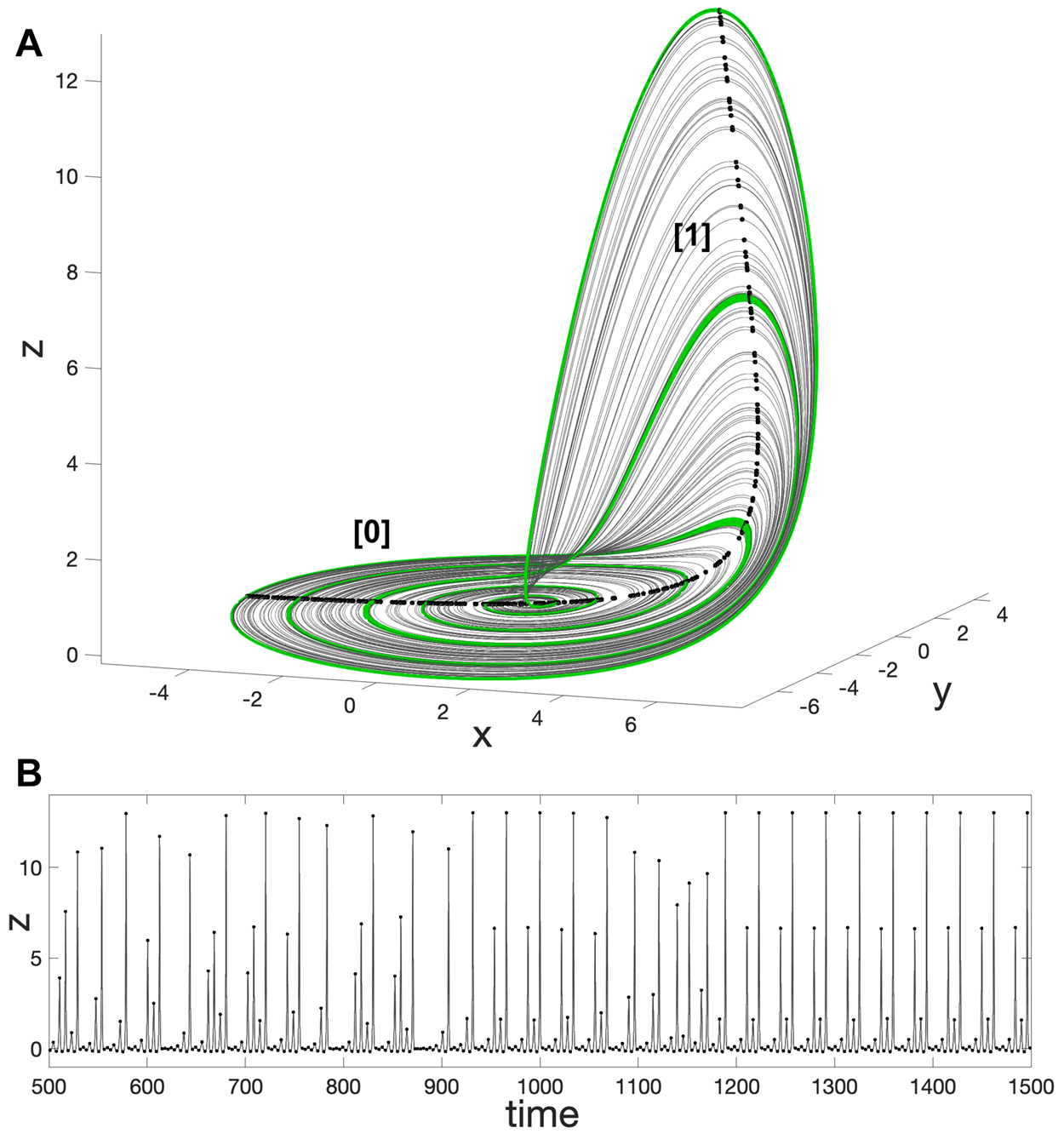}
	\caption{(A) Long chaotic transient (grey) towards a stable periodic orbit (green) in the 3D phase space of the R\"ossler model at $a=0.341$ and $c=4.8$. Black dots indicate the location of $z'=0$ events to generate binary sequences $\{...00011.,.\}$ depending on where the critical events occur below or above some $z$-threshold.  (B) Spiking $z$-variable plotted against time becomes regularized to produce a periodic pattern of low complexity.}
	\label{fig:Rossler}
\end{figure}

The R\"ossler model~\cite{Rossler1976,Rossler1979} is another classical example of deterministic chaos occurring in many low-order systems.  We use the following representation of the R\"ossler model:
\begin{equation}
\dot x = -y - z, \quad \dot y = x + ay, \quad  \dot z = bx +z (x-c),
\label{eq:Rossler}
\end{equation}
 with $x,y$, $z$ being the phase variables, and $a>0$ being bifurcation parameters; here we keep  $c=4.8$ and $b=0.3$ fixed. The convenience of the representation~(\ref{eq:Rossler}) is that  one equilibrium (EQ) state $O_1$ is always located at the origin $(0,0,0)$, while the coordinates of the  second one $O_2$ are given by $(c-ab, b - c/a, -(b - c/a))$.  
 
 The best known feature of the R\"ossler model is the onset of chaotic dynamics due to a Shilnikov saddle-focus \cite{Shilnikov1965,Shilnikov_heritage,Rossler2020}, see Fig.~\ref{fig:Rossler}. It begins with a super-critical bifurcation of the stable equilibrium at the origin, followed by a period-doubling bifurcation cascade as $a$ increases.

Recall that a 3D dissipative system with the Shilnikov saddle-focus cannot produce a {\em genuinely} chaotic attractor but produces a {\em quasi-attractor} \cite{AfrShil1983,GTS93} instead.  Homoclinic tangencies (literately stirred by two saddle-foci in the R\"ossler model) inside such a quasi-attractor cause the emergence of stable periodic orbits in the phase-space through saddle-node bifurcations followed by period-doubling ones. Note that the chaotic attractor in the Lorenz model is categorized as a pseudo-hyperbolic one according to \cite{TS98,TS08,GST97,GKT21}; at larger $r$-values it becomes a quasi-attractor \cite{BykovShilnikov92} because of the presence of various stability islands seen in Fig.~\ref{Lorenz-lyap.fig}. The same is true for the R\"ossler model as illustrated by Fig.~\ref{Rossler-lyap.fig}.

This observation lets us introduce the partition using critical events when the $z$-variable reaches its maximal values on the attractor. The binary sequence $\{k_n\}$ representing a trajectory is computed as follows:
\begin{equation}\label{kneadingsMain}
k_n = \begin{dcases*}
1, & if~~ $z_{max}$~~~~~ > $z_{\rm th}$, \\
0, & if~~ $z_{min/max}$ $\leq$ $z_{\rm th}$, 
\end{dcases*}, 
\end{equation}
where the $z$-threshold can be set relative to the location of the secondary equilibrium state $O_2$: 
$z_{\rm th} = 0.1\,(c - ab)/a$ as was done in Ref.\cite{Rossler2020} for example, or even set it fixed $z_{\rm th} = 0.03$ or $z_{\rm th} = 1.0$ as done in this study. The choice of partition is motivated by its simplicity and may differ if it satisfies the purpose, namely, to distinguish local dynamics concentrated around the Shilnikov saddle-focus at the origin from  large global pathways associated with possible homoclinic excursions in such spiral attractors, see Fig.~\ref{fig:Rossler}    
   With this simple algorithm based on a threshold level $z_{\rm th}$, we can convert the maximal values of the $z$-variable into the binary framework:  using $0$ when a solution of the model turns around the saddle-focus, and  $1$ when it transitions towards the other equilibrium state $O_2$ and back to the origin. Figure~\ref{fig:Rossler}A illustrates the concept in the 3D phase space where the black dots indicate the critical events on the attractors, while Fig.~\ref{fig:Rossler}B represents the time-progression of the $z$-phase variable.

\section{Methods: Complexity measures}
There are several approaches available to assess such a degree. Perhaps the easiest implementation is to convert a time-progression into a binary sequence, apply a compression application, such as gzip, and then to compare the lengths of the original, $L_\text{orig}$, and  compressed, $L_\text{comp}$, files. The compressibility measure can be introduced as $R=1-L_\text{comp}/ L_\text{orig}$ and it is related to the redundancy of information contained in the sequence. As two benchmarks one can choose: (i) non-redundant random Bernoulli sequence with the least compressibility, and (ii) a redundant periodic sequence with the maximal compressibility. A sequence generated by the deterministic chaotic system lies in between of these two benchmark limits. 

In what follows we will employ four different measures to examine the degree of chaos in the Lorenz and R\"ossler models and to determine the stability windows as the parameter is swept within some ranges of interest for these systems. We state from the very beginning, that period-doubling bifurcations are beyond the scope of this examination, and as such our partition designs are not meant to detect such transitions. This can be obviously refined and resolved with additional constrains on the chosen partitions.    
 
 The proposed chaos measures are the source entropy (SE), the Lempel-Ziv complexity (LZ), the Markov transition matrix, and the largest Lyapunov exponent (LE). As the Reader will see, all of these measures work quite well to detect the stability islands in the chaos sea. The key to understanding why this is the case is rooted in the fact that the Lorenz attractor is pseudo-hyperbolic and then all four measures change abruptly with parameter variations, except for stability intervals where they remain constant. Specifically, SE,  LE, and LZ vanish after converging to exponentially stable, and therefore structurally stable, periodic orbits.   
  
\subsection{Block-entropies and the source entropy of symbolic sequences}
A generic measure of complexity of a symbolic sequence is provided by the entropy of the source \cite{shannon1948mathematical}. Given a binary sequence $S$, the Shannon's entropy of words of $m$-symbols long, the so-called block entropy, is defined as \cite{ebeling1992word,ebeling1997prediction}
\begin{equation}
	H_m = -\sum_{\{s_m\}}
	P(s_m)\log P(s_m),
	\label{block_entr.eq} 
\end{equation}
where $P(s_m)$ is the probability of occurrence of a word of length $m$ within the sequence $\S$, and summation is carried over all words of length $m$ occurring with nonzero probability.   The $m$-block entropy, $H_m,$ is interpreted as average information contained in a word of length $m$.
The conditional entropies are defined as in reference\cite{ebeling1997prediction},
\begin{equation}
	h_m = H_{m+1} - H_m, \quad h_0 := H_1,
	\label{cond_entr.eq}
\end{equation}
and provide the average information required for prediction of $(m+1)$ symbol, given that the preceding $m$ symbols are known. The limit of $m\to\infty$ gives the quantity of interest, entropy of the source, 
\begin{equation}
	h=\lim_{m\to\infty} h_m = \lim_{m\to\infty} \cfrac{H_m}{m}.
	\label{sourceent.eq}
\end{equation}

For a dynamical system, the source entropy depends on a particular partition of the phase space, that is, on a rule that maps the evolution of the dynamical system to a discrete symbolic sequence. The Metric or Kolmogorov-Sinai entropy is an upper bound of the source entropy over all finite partitions; the entropies match for the so-called generating partition \cite{eckmann1985ergodic,Sinai:2009}. 

For a symbolic sequence generated by a periodic source, such as a stable periodic orbit  of the Lorenz model, the source entropy equates to 0, while the conditional entropy, $h_m$, drops to 0 when the word length reaches the period, $m=p$, and so $h_p=h_{p+1}=...=0$.
That is, it is enough to observe a periodic sequence for just a period to predict with certainty the next symbol. Equivalently, no new information is gained after one period of periodic sequence is observed. 

For a random or chaotic sequence the source entropy is positive \cite{ebeling1992word,basios2011symbolic}, reflecting the sole fact of uncertainty in prediction of next symbol in chaotic sequence even if the entire prehistory is known. The decay of the conditional entropy, $h_m$, to its asymptotic value $h$ is a generic measure of correlations in the sequence \cite{ebeling1992word,herzel1994finite}. In particular, for a Markov sequence of memory $p$ the conditional entropy converges to the source entropy after exactly $p$ steps, i.e.
$h_{p}=h_{p+1}=...=h$.

The number of words grows quickly with the word length. For example, in the case of Bernoulli binary sequence the number of possible words of length $m$ is $M_m=2^m$. In a numerical experiment the length of a sequence generated by the dynamical system is always finite, which creates a well-known problem in estimation of long-words probabilities \cite{grassberger1988finite,herzel1994finite}. That is, some words cannot be observed or encountered only few times, simply because the sequence is not long enough. In result, the block entropies, $H_m$, are systematically underestimated and finite-size correction must be applied.
As we are interested in an indicator of chaotic or periodic dynamics, in the following we limit the word length to $m=6$ and use $h_6$ for estimation of the source entropy. We collect long sequences, $N\gg2^6$, and use the finite size correction \cite{herzel1994finite},
\begin{equation}
	H^\text{observed}_m \approx H_m - \frac{M_m-1}{2N},
	\label{entr_obs.eq}
\end{equation}
where $	H^\text{observed}_m$ is the observed $m$-block entropy calculated from the observed  sequence of length $N$, $H_m$ is the true entropy, $M_m$ is the number of distinct $m$-words, $M_m \le 2^m$. 

\subsection{The Lempel-Ziv complexity}
Given a binary sequence, the Lempel-Ziv  complexity (LZ) is related to the number of substrings of increasing lengths that the given sequence is made of.  For example, for the sequence, $\{1|0|10|010|0101|11|110\}$, scanned from left to right, that number $s$ is 8.
Then $LZ$-complexity can be defined as 
\begin{equation}
	LZ=s \log(N)/N,
\end{equation}
where $N$ is the length of the sequence.  Indeed, the LZ algorithm can be used as an effective estimator of the source entropy \cite{ziv1978compression}. Thus, we expect the source entropy estimate, $h_6$, and the LZ  to be strongly correlated, as the parameter of a dynamical system varies.

\subsection{The Markov transition matrix}
The simplest measure is introduced when it is based on elements of the transition matrix of the single-step Markov process underlying the symbolic sequence,
$$\textbf{M}= \left [ \begin{matrix} p_{11} & p_{01} \\  p_{10} & p_{00}
\end{matrix} \right ],$$
where $p_{11}$ and $p_{10}$ are probabilities of flop and flip, respectively, i.e., of transitions $1 \to 1$ and $1 \to 0$. Because the Lorenz attractor is symmetric,  $p_{11} \simeq p_{00}$, which implies that  $p_{10} \simeq  p_{01}$ as well. The situation is different when the  attractor is a stable asymmetric periodic orbit, or either one in a pair of asymmetric chaotic attractors that emerged through a period-doubling cascade. The transition probabilities $p_{11}$ and $p_{01}$ can be compared with $0.5$, e.g. with Bernoulli trial of an unbiased coin. For example, the Markov matrices for the symmetric stable periodic orbits shown in Figs.~\ref{fig2} and \ref{fig3} are, resp., the following 
 $$  \left [ \begin{matrix} 4/5 & 1/5 \\  1/5 & 4/5
\end{matrix} \right ]  \quad \mbox{and} \quad \left [ \begin{matrix} 3/4 & 1/4 \\  1/4 & 3/4
\end{matrix} \right ], 
$$
and will remain such within the corresponding stability windows in the parameter space. 
  
\subsection{The largest Lyapunov exponent}
  \begin{figure}[t!]
	\includegraphics[width=0.9\linewidth]{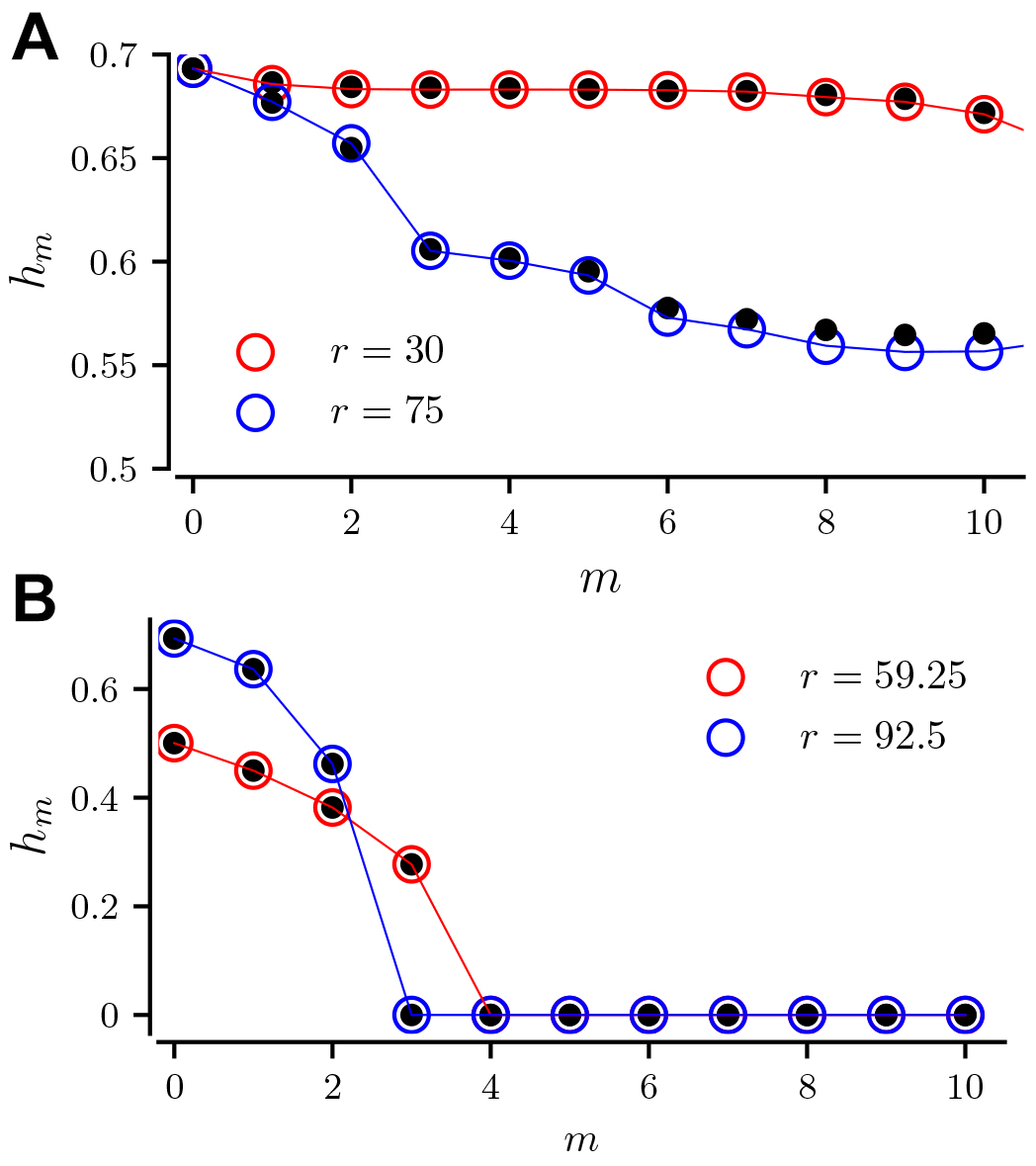}
	\caption{Conditional entropy versus the words length for the Lorenz model.
		(A) Conditional entropy for chaotic sequences with  $N=16988$ for $r=30$, and $N=28475$ for $r=75$. (B) Conditional entropy for periodic sequences with length $N=25672$ for $r=59.25$,  and $N=32058$ for $r=92.5$. Open circles refers to the conditional entropy $h_m$ estimated from a sequence generated by the genuine Lorenz model. Filled black circles shows $h_m$ estimated from sequences generated by the trained reservoir computer whith the same lengths of binary sequences.}
	\label{entro1.fig}
\end{figure}

A Lyapunov exponent is meant to indicate how quickly nearby trajectories may converge/diverge in the phase space.  The sum of positive Lyapunov exponents is related to the Kolmogorov-Sinai (KS) entropy via Ruelle's inequality, $\text{KS} \le \sum_{\Lambda_i>0} \Lambda_i$ \cite{ruelle1978inequality}. Since the Lorenz system is a strongly dissipative system that may have a single positive Lyapunov exponent on the chaotic attractor, then the largest Lyapunov exponent (LE) is directly related to the KS entropy,  $\text{KS} \le  \Lambda$.
Thus, the measures introduced for symbolic dynamics can be compared against the LE,
\begin{equation}
	\Lambda = \lim_{t\to\infty}\frac{1}{t}\log \frac{|\delta \textbf{x}(t)|}{|\delta \textbf{x}_0|}, \nonumber
\end{equation}
that measures average rate of convergence or divergence between two trajectories. So, $\Lambda<0$ means a trajectory converges to a stable equilibrium state; $\Lambda=0$  corresponds to the case of a stable periodic orbit along which the distance between two solutions does not change over its period. The case $\Lambda>0$ indicate that solutions of the system under consideration run on some chaotic attractor.

For a fair comparison with symbolic-sequence measures, such as the source entropy, the LE should be normalized using a characteristic time of the system,  $\tau$,  so that the quantity $\lambda = \Lambda \tau$ becomes dimensionless. For the Lorenz attractor case, as the characteristic time we use the mean dwelling time intervals between the events $\dot x=0$, such as ones marked as the red dots on the butterfly wings in Fig.~\ref{fig1}. One can observe from Eqs.~(\ref{eq:Lorenz}) that increasing the parameter $r$ speeds up the time derivative $\dot y$ and hence $\dot x$, which results in shortening of the dwelling time intervals between transition events which is compensated by the growing size of the Lorenz attractor, compare the coordinate axis in Figs.~\ref{fig1} and \ref{fig:L75}, for example. Similarly, for the R\"ossler model below we used the dwelling times between similar critical events.  

\begin{figure*}[t!]
	\includegraphics[width=0.9\linewidth]{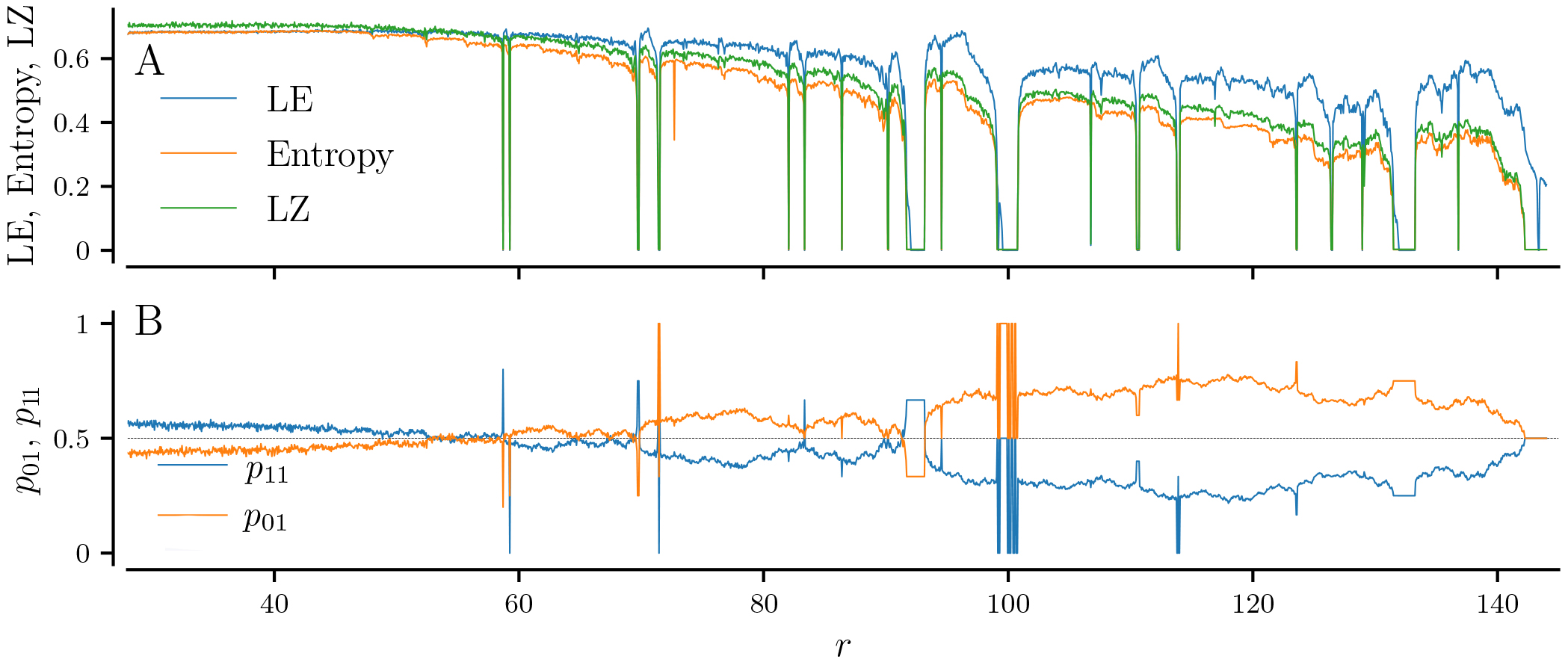}
	\caption{Complexity measures vs the parameter $r$ of the Lorenz model.
		For each parameter $r$ value, symbolic sequences of the length $10^4$ were collected.
		A: The largest Lyapunov exponent (LE, blue), the source entropy estimate, $h_6$ (Entropy, orange), and  Lempel-Ziv complexity (LZ, green) vs the parameter $r$. The largest Lyapunov exponent was normalized to the average inter-symbol interval, $\lambda T$.
		B: Markov probabilities $p_{11}$ and $p_{01}$ show sharp peeks and plateaus  that are also indicative, resp., of narrow and wide stability windows with periodic orbits within. 
	}
	\label{Lorenz-lyap.fig}
\end{figure*}

\section{Data and Results}
 
Figure ~\ref{entro1.fig} presents the conditional entropies for chaotic and periodic regimes of the Lorenz model.  The first examples in Fig.~\ref{entro1.fig}A illustrate the practicality of the conditional entropy in distinguishing two chaotic regimes: larger  $h_m$ value for the aperiodic binary sequence generated at $r=30$ compared to that at $r=75$ is a direct indicator of the higher degree of uncertainty of symbols predictions in the former. This implies that one can treat chaos at nearly canonical parameter value is superior than that at $r=75$, see Fig ~\ref{fig:L75}A. One indirect justification for that is the Lorenz model at  $r \ge 31$ no longer exhibit genuinely chaotic attractors but quasi-attractors with homoclinic tangencies causing the emergence of stable periodic orbits such as one occurring at $r=68.75$ and shown in Fig.~\ref{fig2}. This should imply that binary sequences generated by chaotic solutions of the Lorenz model at this parameter range include multiple ``laminar'' or periodic substrings leading to lower values of conditional entropy.  We will proceed with more direct arguments in the next section.

As for stable periodic orbits generating stable periodic orbits after some chaotic transients, the conditional entropy quickly vanishes when the word length reaches the period, see Fig.~\ref{entro1.fig}B. The distinction between two periodic regimes is clearly captured by the entropy. First, the obvious observation is that the stable periodic orbit observed at $r=59.25$ has a longer period than the orbit found at $r=92.5$. A less obvious observation is related to the value of $H_1$ (or $h_0$), that is the entropy of words of length 1. 
For $r=92.5$ $H_1$ is given by $\log 2$, reflecting that probabilities of "0" and "1" are the same. This indeed reflects the symmetry of underlying stable periodic orbit. On the contrary, for $r=59.25$, the entropy $H_1$ is significantly smaller than $\log 2$, reflecting the asymmetry of the underlying  periodic orbit.

\begin{figure*}[t!]
	\includegraphics[width=0.9\linewidth]{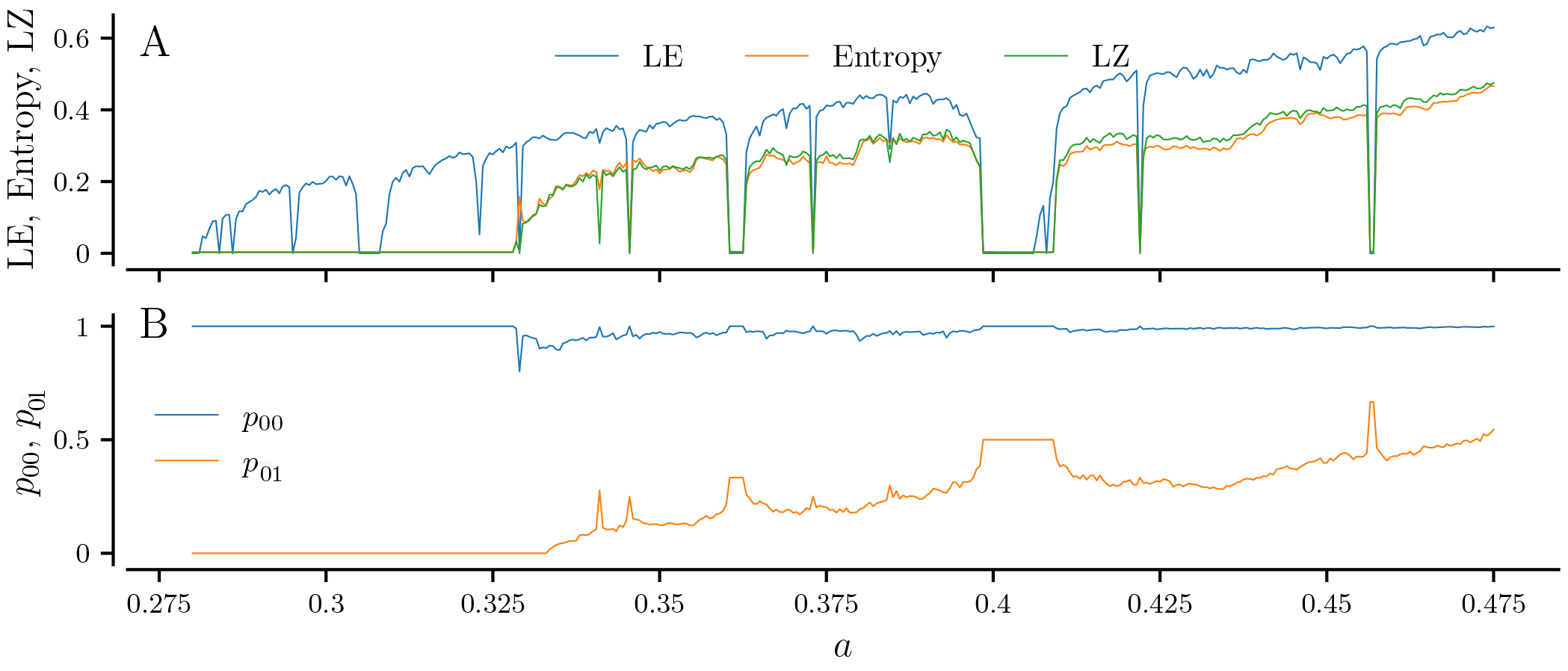}
	\caption{Complexity measures vs the parameter $a$ of the R\"ossler model.
		For each parameter $a$ value, symbolic sequences of the length $10^4$ were collected.
		A: The largest Lyapunov exponent (LE, blue), the source entropy estimate, $h_6$ (Entropy, orange), and  Lempel-Ziv complexity (LZ, green) vs the parameter $a$. The largest Lyapunov exponent was normalized to the average inter-symbol interval, $\lambda T$.
		B: Markov probabilities $p_{00}$ and $p_{01}$. Sharp dropping peeks and plateaus  are indicative of narrow and wide stability windows with periodic orbits within. 
	}
	\label{Rossler-lyap.fig}
\end{figure*}

Figures~\ref{Lorenz-lyap.fig} and \ref{Rossler-lyap.fig} with several charts are the graphical culmination  of our simulations in the first part of this study dealing with complexity measure of  dynamics demonstrated by  the Lorenz and R\"ossler models. As we said above that leave the chaos due to period-doubling bifurcations apart from our consideration, as follows directly from the choice of partitions in either case.   

So, Figure~\ref{Lorenz-lyap.fig} represents a $r$-parametric sweep of dynamics of the Lorenz model starting with the classical parameter value $r=28$: its panel~B depicts empirically the way the complexity measures vary and all correlate with $r$-parameter increase.  The low panel B in Fig.~\ref{Lorenz-lyap.fig} shows variations of two (asymmetric) entries, $p_{11}$ and $p_{01}$, within the same $r$-range. Note that here $p_{11}$ and $p_{01}$ taken from different columns of the transition matrix and hence do not add up to $1$ always, are used to detect stability windows with paired mirror-asymmetric periodic orbits within.    

Let us first observe the sharp peeks and plaques in these sweeps where all three measures, LE, SE, and LZ  drop down to zero. One can conclude that such $r$-intervals correspond to narrow or wider stability windows within which the Lorenz model will eventually demonstrate periodic dynamics. The corresponding periodic orbit can be symmetric like those shown in Figs.~\ref{fig2} and \ref{fig3}. By analyzing the values of $p_{11}$ and $p_{01}$-probabilities from Fig.~\ref{Lorenz-lyap.fig}B, one can deduct,  without visual inspection of Fig.~\ref{fig3} that the plateau around $r=92.5$ correspond to the symmetric stable orbit, while sharp peeks near $r=100$ are indicative of the stability windows with a pair of co-existing stable orbits within, to either one the transient converges eventually.  

One can see from Fig.~\ref{Lorenz-lyap.fig}A that at the initial stage prior to two first stability windows around $r=59$, the values of the SE and LE agree quite well with each other, while after the LE-curve starts diverging from the LZ and SE curves.  As we pointed out above, increasing $r$-parameter increases on average
the time-derivatives of the $x$ and $y$-phase variables, accelerating the time course in the Lorenz model. 
While we used the average time interval between the critical events at assigning binary symbols as the normalizing factor for comparison of LE and SE,
the distinction is still noticeably growing with $r$-parameter increase.
A possible reason can be that at larger $r$-values the dynamics of the Lorenz model is no longer flip-flop due to homoclinic bifurcations of the saddle at the origin but becomes mostly determined by period doubling bifurcations that our binary framework is not designed to account.
In other words, the used binary partition is far from a generating partition. Such regimes will have inclusions of longer, self-similarly regularized blocks like $....111000111000...$ that result in lower entropy- and LZ-based measures, whereas the average rate of divergence of two nearby trajectories in the phase space remains basically intact. 

\begin{figure*}[htb!]
	\includegraphics[width=0.999\linewidth]{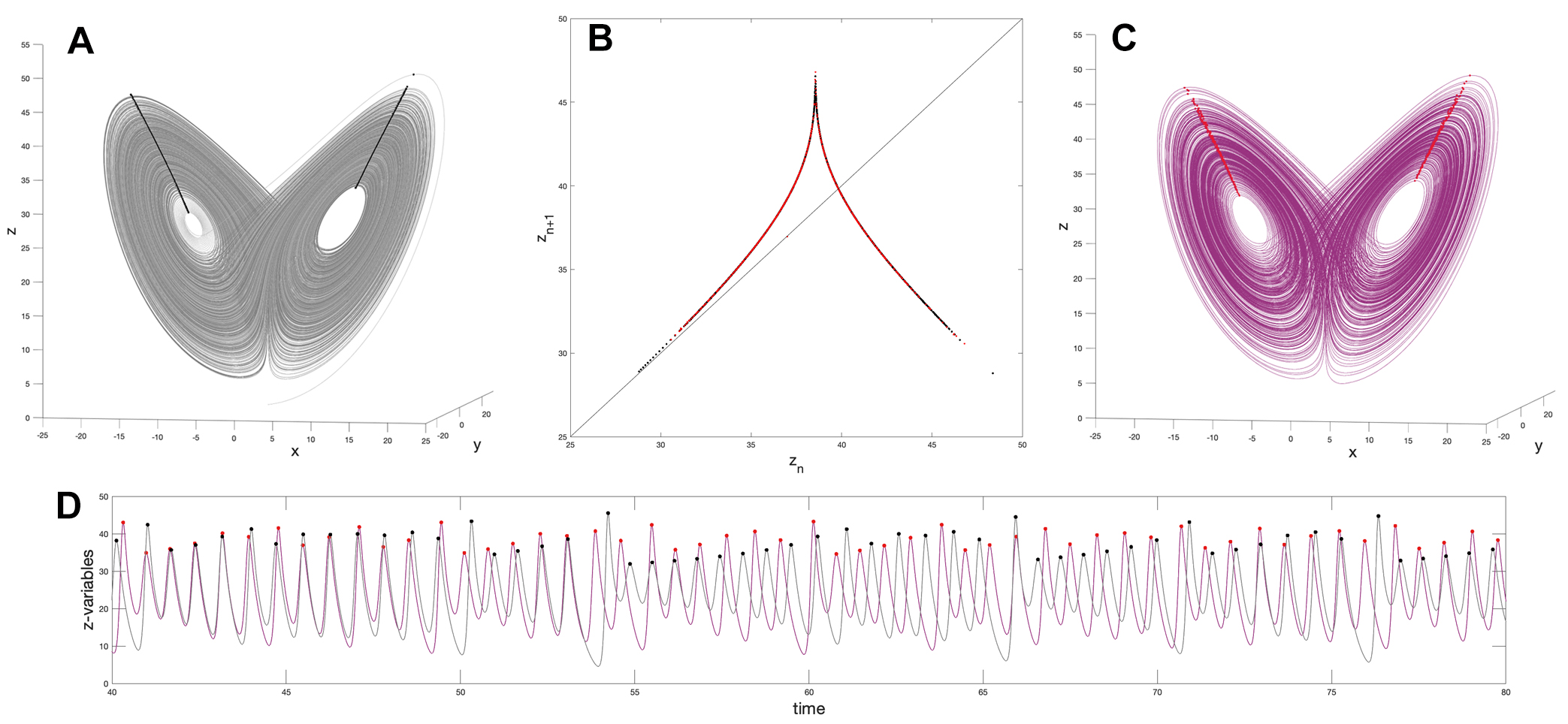}
	\caption{(A) The Lorenz attractor at  $r=28$. Superimposed black dots are the $z$-critical events well-aligned on two line segments transverse to the butterfly wings in the phase space. (B) The $z$-variable plotted against time is superimposed with local maxima given by $\{z'=0|\, z''<0\}$ to generate the Poincar\'e return maps shown in (C). (C) 1D maps $T: z_n \to z_{n+1} $ generated by sequential pairs of local maxima of the $z$-variable in the Lorenz model (black dots), while the map filled out with red dots are generated from the ML-emulated trace (red) shown in (D). The Lorenz $z$-map has a characteristic cusp-shape (due to high outbursts of the separatrix of the saddle, see (A)) and an emerged fold at $r=75$ giving rise to stable orbits due to tangent (SN) bifurcations. The emulated map adequately captures some characteristics of the original on average. }
	\label{fig:L28}
\end{figure*}

 Figure~\ref{Lorenz-lyap.fig}A  suggests non-ambiguously that chaos in the Lorenz model at low $r$-values is more unpredictable and homogenous due to shorter flip-flop dynamics caused by homoclinic bifurcations of the saddles than at higher parameter values where the dynamics is mostly dominated by longer recurrent patterns due to period-doubling bifurcations. Figure~\ref{Lorenz-lyap.fig}B supports this assertion that the complexity is maximized till the entries of the transition matrix remains close to 1/2 on the left in the sweep before the very first stability islands.

One-parametric sweeps of the  largest Lyapunov exponent, the LZ-complexity, and the source entropy for the R\"ossler model are presented  in Fig.~\ref{Rossler-lyap.fig}A, while Fig.~\ref{Rossler-lyap.fig}B depicts the way the Markov transition probabilities change accordingly with variations of the $a$-parameter. One can see from its initial phase $a<0.33$ that complexity measures stay close to zero, while the LE is indicative of the origin of chaotic dynamics alternating with periodic dynamics in the phase space of the R\"ossler model.  
The examination of Fig.~\ref{Rossler-lyap.fig}A reveals that the probability $p_{00}=1$ till $a<0.33$, which indicates by the phase space partition that the dynamics of the R\"ossler model remains ``flat'' near the origin (see Fig.~\ref{fig:Rossler}) that became a saddle-focus after a super-critical Andronov-Hopf bifurcation followed by period-doubling cascades through which the dynamics becomes weekly chaotic as indicated by the positive largest Lyapunov exponents.

For $a>0.33$, the Shilnikov saddle-focus starts contributing to more developing chaotic dynamics  when its 2D stable manifold bends so that trajectories come closer to the proximity of the origin. With large $a$-values the chaotic attractor appears to look more like iconic spiral chaos in the R\"ossler model.  
The growing values of the  complexity measures along with the $p_{01}$-probability in the transition matrix approaching 1/2 from below are all synchronous indicators that the degree of chaos and unpredictability on the R\"ossler model increases as it transitions further from period-doubling type towards spiral and multi-funnel chaos with countably many homoclinic orbits originated by and stirred by both Shilnikov saddle-foci in the model, see Ref.\cite{Rossler2020} for more details. As was said earlier, in addition to chaos, the Shilnikov homoclinic bifurcations lead to the emergence of homoclinic tangencies, and hence the abundance of  saddle-node bifurcations that determine the borders of multiple stability windows in the parameter space that are populated by sequential period-doubling cascades within.  This is well-seen in the sweep in Fig.~\ref{Rossler-lyap.fig}. The reader is welcome to consult with Refs.~\cite{xing2014symbolic,pusuluri2018homoclinic, pusulurihomoclinicCNSNS,XPSh2010} that elucidate computationally the contribution of the Shilnikov saddle-focus bifurcation in the formation of the so-called T-points for heteroclinic connections between saddles and saddle-foci in various Lorenz-like systems and the role of inclination-flip bifurcations in such unfolding that amplify homoclinic tangencies in producing Shilnikov flames -- stability islands nearby such T-points.

To conclude this section we re-iterate that the conditional entropy captures well the parameter dependence of the dynamics, showing qualitatively the same dependence as the largest Lyapunov exponent, and LZ-complexity as de-facto depicted by the parametric sweeps of the Lorenz and R\"ossler model. Let us remind that for a fair comparison in the largest Lyapunov exponent was normalized by the average time intervals between the critical events used in the binary framework design. Finally, the transition probabilities cab be effectively used as indicators to detect stability  windows in parametric sweeps as well.

\section{Reservoir-computing}

\begin{figure*}[htb!]
	\includegraphics[width=0.999\linewidth]{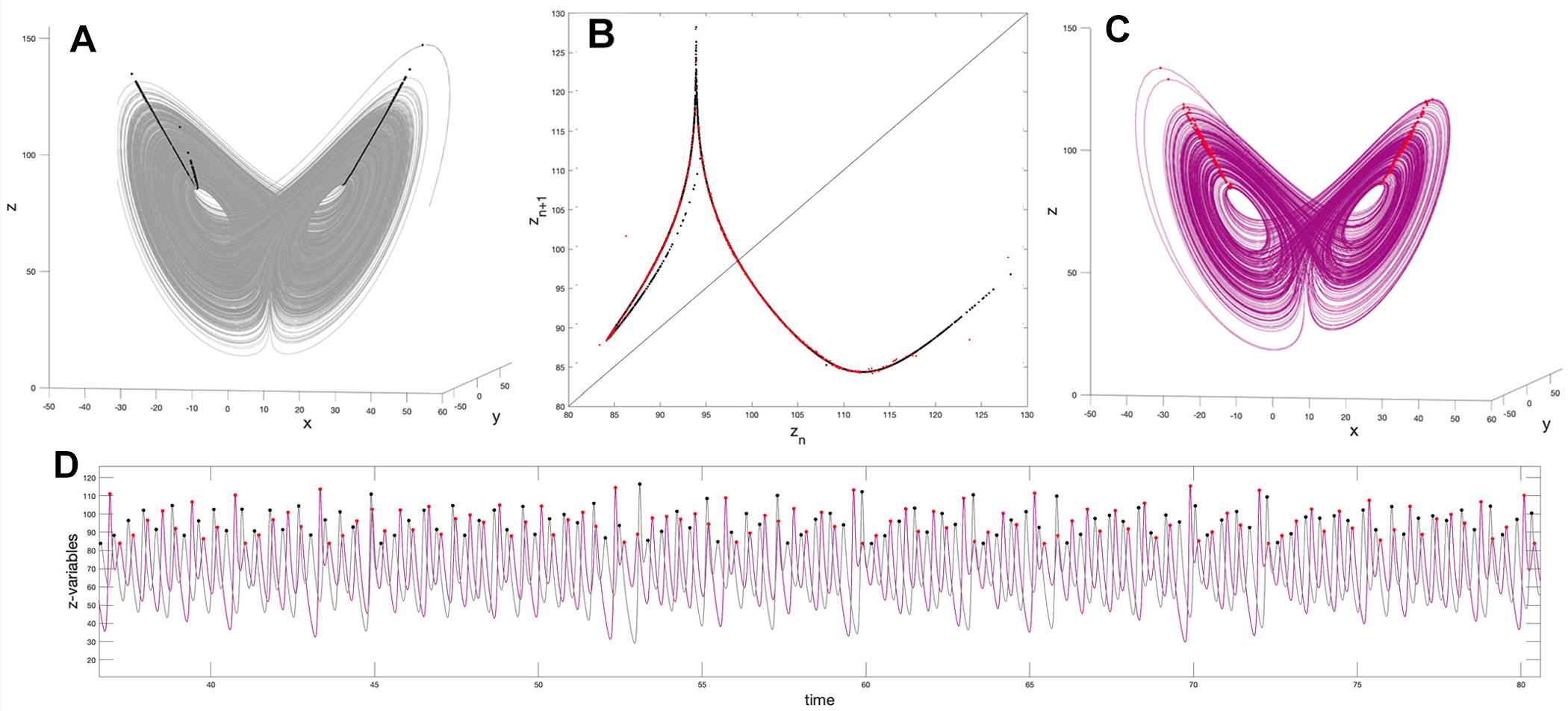}
		\caption{(A) Chaos in the Lorenz model is shown at $r=75$. Superimposed black dots are $z$-variable critical events form two 1D hooks on the bending butterfly wings in the phase space. (B) The $z$-variable plotted against time is superimposed with local maxima given by $\{z'=0|\, z''<0\}$ to generate the Poincar\'e return maps shown in (C). (C) 1D maps $T: z_n \to z_{n+1} $ generated by sequential pairs of local maxima of the $z$-variable in the Lorenz model (black dots) and one (red dots) from the ML-emulated trace (red) shown in (D). The Lorenz $z$-map (no longer 1D due to loss of strong contraction) has a characteristic cusp-shape (due to high outburst of the separatrix of the saddle (see (A)) and an emerging fold at $r=75$ to give rise to stable orbits due to tangent (saddle-node) bifurcations. The emulated map captures adequately some characteristics of the original on average. } 
	\label{fig:L75}
\end{figure*}

Symbolic dynamics can be used to study the properties of systems identification techniques for chaotic systems. Echo State Networks (ESNs), a form of reservoir computing (RC), are one such method for systems identification that has recently become popular in the nonlinear dynamics community \cite{Jaeger01,reservoirreview}. In the following section we apply a fidelity test to  echo state networks, trained for time series prediction. The fidelity test includes two components: a quantitative test through complexity measures, specifically conditional entropies of various block sizes LZ complexity, and a qualitative test through canonical 1D Poincar\'e return maps :  $T$: $z_{max}^n \to z_{max}^{n+1}$.\\

In an ESN, the input is embedded into the reservoir’s vector space by a random linear transformation. The resulting input and reservoir vectors are summed and passed through some generic point-wise nonlinearity, hyperbolic tangent in this case. This is derived from the discretization of a nonlinear leaky integrator, yielding the update relation for the reservoir \cite{Jaeger01}:

\begin{equation}
	R_{n+1} = (1-\alpha)R_n  + \alpha \tanh (W_RR_n+W_{in}y_n),
\end{equation}
where $R_n$ is the reservoir and $y_n$ is either the teacher input, or the last prediction: $W_{out}R_n$. 

The input reservoir and input matrices are $W_R$ and $W_{in}$, respectively. The final prediction is obtained from $R_n$ as
\begin{equation}
	x_{n} = W_{out}R_n,
\end{equation}
where $W_{out}$ is the matrix representation of the best fitting linear transformation of the reservoir sequence, found by ridge regression. Detailed methods for both Lorenz and R\"ossler models can be found in the following sections. The results are summarized for each and followed by a brief discussion below.

\subsection{Lorenz model's replica}

The results of the quantitative test are summarized in  Fig.~\ref{entro1.fig}. Shown in circles are the values of the conditional entropy $h_m$ estimated from binary sequences generated by the Lorenz model, while the  black dots represents the $h_m$-values estimated from binary sequences generated by the trained reservoir computer.

The results of the qualitative test are illustrated in the sequence of two Figs.~\ref{fig:L28}, and \ref{fig:L75}. Each figure includes four panels: two panels on the left and right depict the phase space projections of the original Lorenz attractor (grey) and its RC reconstruction (purple). The black/red dots in these panels indicate the critical events when the $z$-variables become maximized. The bottom panel show the time progression of the original $z$-variable overlapped with that of the surrogate, as well as the $z_{max}$ map in the middle panel. This is a standard approach that is also applicable to RC surrogate systems ~\cite{pathak2018model,pathak2017using}.

Let us begin with the return maps shown Fig.~\ref{fig:L28}B. One can see that both maps, ordinal and surrogate, produce the expected cusp-shape at $r=28$ with the slope being steeper then $\pm 1$ by comparing the graph with the bisectrix or the $45^\circ$-line presented in this figure as well, i.e., the map is an expansion.  One can  see that the cusp graph populated by the $z{max}$-values (black dots) extracted from the Lorenz model is exceptionally well fitted by those (red dots) from the RC surrogate.    

Such a cusp map constitutes a de-facto computational proof that the Lorenz attractor has no stable orbits and is made only of countable many unstable/saddle periodic orbits whose 2D stable and unstable manifolds cross transversely in the 3D phase space of the model. This assertion follows also indirectly from the visual inspection of the dots, given by the condition $x'_{max/min}=0$, that fill out two straight line segments on the opposite sides of the symmetric wings of the Lorenz attractor shown in Fig.~\ref{fig:L28}A while its RC surrogate is shown in Fig.~\ref{fig:L28}C. 
  
  \begin{figure*}[htb!]
	\includegraphics[width=0.999\linewidth]{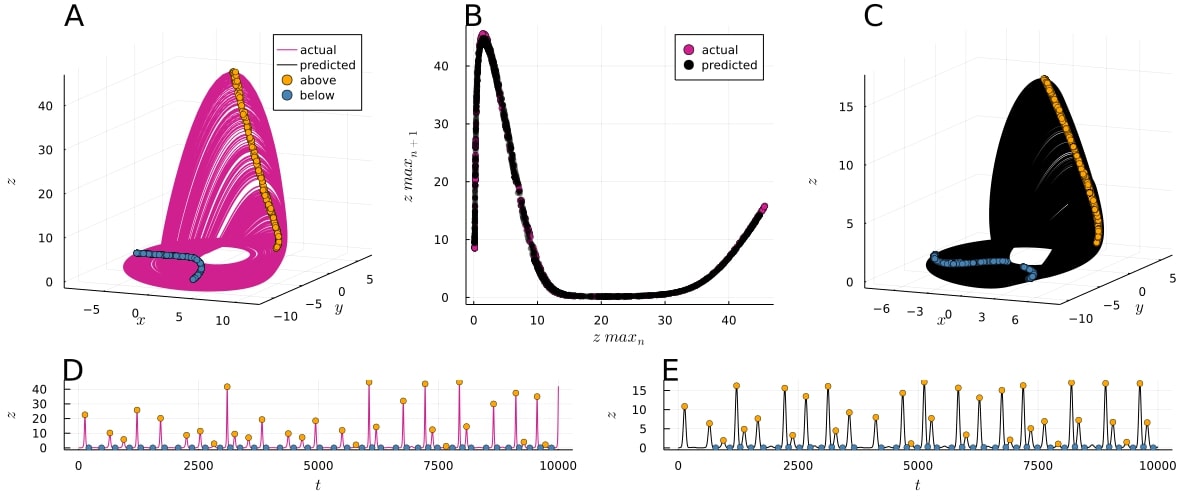}
		\caption{True solution (A) and predicted solutions (C) for R\"ossler attractor are shown with  orange dots demarcating local maxima with $z>1$ and blue dots demarcating local minima. Accompanying time series are given in (D) and (E) respectively. (B) shows the characteristic $z_{max}$-return maps.}
	\label{fig:rosslerfig}
\end{figure*}
\begin{figure}[h!]
	\includegraphics[width=0.9\linewidth]{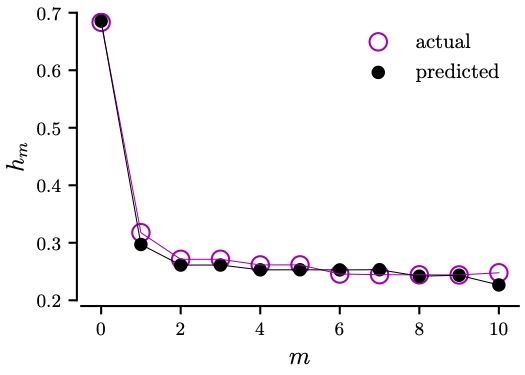}
	\caption{Conditional entropy versus the words length for the R\"ossler system and its reservoir computer clone. The parameters are the same as in the previous figure.
		Open circles refers to the conditional entropy $h_m$ estimated from a sequence generated by the genuine R\"ossler system. Filled black circles shows $h_m$ estimated from sequences generated by the trained reservoir computer with the same lengths of binary sequences.}
	\label{fig:rossler_entro.fig}
\end{figure}

This is no longer the case at greater $r$-values in the Lorenz model, as was first shown computationally in Ref.~\cite{BykovShilnikov92}, namely, stable and unstable manifolds of saddle periodic orbits populating such quasi-attractors may no longer cross transversely, thereby causing homoclinic tangencies that give rise to the onset of stable periodic orbits (inside stability windows in the parameter space) through saddle-node bifurcations. 

Let us first discuss the change that the shape of the return map undergoes at higher $r$-values, specifically the maps at $r=75$ are contrasted with corresponding maps with $r=28$.  One can observe first of all from Fig.~\ref{fig:L75}B  that the map becomes extended with the bending section on the right. This fold is where the map looses its expansion property due to the flat portion with zero derivative. This is a direct cause for forthcoming saddle-node bifurcations producing stable periodic orbits which can coexist with the chaotic pseudo-hyperbolic subset (comprised of countably many periodic orbits) in bistability or become global attractors in the phase space. The images of such [close to, due to strong contraction in the transverse direction] 1D folded return map can be also observed in the phase space projections of the chaotic set in the Lorenz model~Fig.~\ref{fig:L75}B and its RC surrogate in Fig.~\ref{fig:L75}C. These maps populated by the black and red dots (given by $x'_{max/min}=\ 0$) demonstrate the distinctive hooks/folds where the wings of the attractors bend in the phase space, unlike that case $r=28$ where they remain straight and flat.

The ground truth trajectory was integrated with Matlab's ODE45 with fixed timesteps at .005. The Lorenz ESN was based on Matlab code by M. Lukosevicius \cite{lukovsevivcius2012practical}. Detailed parameters can be found on our github, cited below \cite{github}.

\subsection{R\"ossler model's replica}
The qualitative results for the R\"ossler model's fidelity test are displayed in Fig. \ref{fig:rosslerfig}. Panels A and D correspond to the real model and panels C and E correspond to the predicted trajectories. The $z_{max}$-return map fits remarkably well, shown in panel B. The geometric distortion of the predicted attractor is an artifact of smoothing which applied to ameliorate the detection of "false" critical points when calculating the symbolic dynamics, but the distortion was included for illustrative purposes. 

Quantitative results can be seen in Fig \ref{fig:rossler_entro.fig}. The conditional entropies show tight agreement, even as block size increases, demonstrating the effectiveness of the appraoch. 

\begin{figure*}[htb!]
	\includegraphics[width=0.999\linewidth]{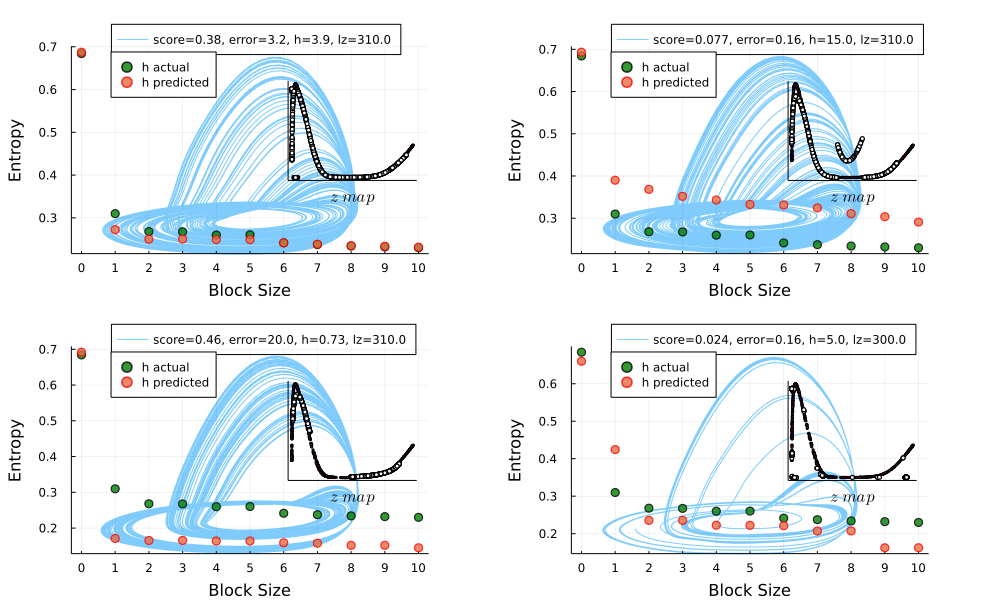}
		\caption{3D projections of four surrogate attractors generated by the R\"ossler random search, each displaying different qualitative properties. The top legend shows the results of the short-run scoring system, and the trajectory in phase space (background) illustrates each. The actual and predicted conditional entropies are plotted, along with the $z_{map}$-return maps superimposed on the right.}
	\label{fig:searchsample}
\end{figure*}

We note that because of slow-fast nature of the R\"ossler model and the lack of symmetry,  it was not easy to achieve such a good representation, compared to with the Lorenz model. Several additional steps were necessary to obtain a close replication of the R\"ossler model by ESN. In the particular example of Fig.\ref{fig:rosslerfig}, 	the true solution of the R\"ossler system was obtained for the parameters $a=0.341$, 
$b =0.3$, $c = 4.8$ and was integrated with an adaptive 3rd-order Bogacki-Shampine method \cite{diffeq,julia}. The solution was stored at time steps of 0.02. Noise sampled from a normal distribution with mean 0 and variance 0.3 was subsequently added. Selecting an adequate noise level proved critical in obtaining consistent results, but setting it too high tended to yield periodic orbits rather than chaotic attractors.

A random search over a broad space of hyper-parameters of the ESN was conducted at a reservoir dimension of 80. Each ESN generated in the search was trained over 30000 time steps after discarding a 2000 time step transient using the Julia's library ReservoirComputing.jl \cite{esnjl}. and the score was computed over a test prediction of an additional 30000 time-steps according to the following system:
\begin{itemize}
\item Score = LZ$\times$Transient$\times$Entropy + Symmetry;
\item Lz is the absolute value of the difference between the LZ complexity of the test and predicted sequences;
\item Transient is the squared error over the first 1000 steps, testing for initial synchronization;
\item Entropy is the sum of the differences in block entropies over block sizes ranging from 2 to 10;
\item The symmetry term is designed to disqualify symmetric systems which were commonly generated. It simply adds $10^9$ if the trajectory dips below z = -1. 
\end{itemize}

Symbolic sequences were calculated as a sequence where local min are marked as “0” and local maxima with $z>1$ are marked as “1”. Thus when the trajectory completes a rotation around the origin without crossing the z=1 plane, a “00” is recorded. Computing the local minima for the predicted trajectories was complicated by the existence of “false” local minima, which were ameliorated by three applications of a moving average filter with window 20. This seemed to outperform other combinations of higher order Savitzky-Golay filters, although no rigorous analysis was performed. Critical points were computed using FindPeaks.jl with minimum prominence set at 0.1 with a minimum distance of 150. 

Each ESN generated for the search was trained over 30000 time steps after discarding a 2000 time step transient using ReservoirComputing.jl and the score was computed over a test prediction of an additional 30000 time-steps. The top 50 scoring trajectories were inspected in phase space by eye alongside a plot of their associated block entropies calculated over a longer trajectory of 168,000 steps. The top candidate was chosen and the $z_{max}$ return map and long run conditional block entropies were computed, with the final entropies calculated over 9,968,000 time steps yielding sequences of lengths 69,994 and 66,322 characters for the true and simulated solutions respectively.


\section{Discussion}

In conclusion we would like to reiterate that the symbolic approach used for the characterization of complex chaotic dynamics, followed by the application of the quantitative measure such as 
Shannon's block entropies

The quick block-entropy method, or a longer Lempel-Ziv complexity approach or even a simple algorithm based on a Markov transition matrix works exceptionally well for the given purpose. 
One may argue that its only ``flaw'' is the choice of the proper partition for the phase space to introduce the symbolic description.  While such a partition for the Lorenz model seems evident when we study the homoclinic portion of its parameter space that results in chaotic flip-flop patterns, nevertheless it becomes ineffective in the period-doubling region of the parameter space for larger $r$-values. We stress again that period-doubling bifurcations were out of the scope of our study focused specifically on homoclinic ones in both system: the Lorenz and R\"ossler models. Lastly, let us add that this approach should be viewed not as a substitution but a novel and welcome addition to the computational approached based on the Lyapunov characteristic exponents. Our findings are meant to validate this assertion. 
  
There are two pivotal impressions gained as a result of our experimentation with ESNs. First, it must be acknowledged that while the echo state networks presented here show excellent statistical agreement, that is not necessarily an inherent property of the ESN framework. Curiously, the Lorenz model showed immediate statistical agreement upon finding a parameter set. The R\"ossler model, on the contrary, did not. Many R\"ossler ESN surrogates had a relatively long initial synchronization, but performed poorly statistically. One potential reason for that is that the R\"ossler model is a slow-fast systems, where the $z$-variable happens to be the fast one with sharp low and high peaks. Furthermore, the Lorenz model was much easier to train then the R\"ossler one. R\"ossler surrogates often displayed inappropriate symmetry, collapsed into periodic orbits of various lengths, or converged to visually similar but qualitatively distinct systems. A few examples with accompanying statistical portraits are displayed in the concluding Fig.~\ref{fig:searchsample} below.

We found ESN's to have poor stability in general, and require a significant amount of tweaking. We also found that the generally accepted guidelines for hyper-parameter selection were not particularly useful. The spectral radius, for example, did not have any obvious relation to performance. The best surrogates were found widely distributed throughout the hyper-parameter search space. That said, we make no claim to expertise in the domain, but do note it requires patience and luck. Perhaps someone else with more skill can do better. As such, a method for systematic reservoir construction with improved stability and intuitive parameters would be a welcome contribution.

The second worthy note is that the fidelity test proved to be highly valuable as a tool for automatically filtering ESN search results. Error alone is not a sufficient test for chaotic systems due to sensitivity to initial conditions. The incorporation of statistical quantification into the search criteria proved to be immensely practical, albeit imperfect. The further development of search techniques based on symbolic dynamics, perhaps combined with cross validation or Monte-Carlo search, is a promising area for future research which would yield immediate practical value for chaos prediction methods of many varieties. 

In any case, we can attest that reservoir computing passed all our fidelity tests, a qualitative one based on the methods of dynamical systems and a quantitative test based on statistical properties of chaotic solutions of the original systems such as the Lorenz and R\"ossler models.


\section*{Acknowledgements}\label{sec:s1}
A.~Shilnikov acknowledges partial funding support from the Laboratory of Dynamical Systems and Applications NRU HSE, of the Ministry of Science and Higher Education of Russian Federation, grant No. 075-15-2019-1931. 

\section*{Data/Code Availability}
All codes and data are freely available upon request.

\section*{References}

%

\end{document}